\documentclass{article}
\usepackage{latexsym, fullpage}

\usepackage{mathrsfs}
\usepackage[english]{babel}
\usepackage[utf8]{inputenc}
\usepackage{graphicx, color}
\usepackage[colorinlistoftodos]{todonotes}
\usepackage{fullpage}
\usepackage{graphics, color}
\usepackage{psfrag}
\usepackage{amsmath}
\usepackage{enumerate}
\usepackage{amsthm}
\usepackage{hyperref}
\usepackage{calrsfs}
\usepackage{tikz}
\usepackage{rotating}
\usepackage{amssymb}
\usepackage{amsmath}

\usepackage{authblk}
\usetikzlibrary{positioning,chains,fit,shapes,calc}

\newtheorem*{dfn}{Definition}

\newtheorem{theorem}{Theorem}[section]

\newtheorem{conj}[theorem]{Conjecture}
\newtheorem{problem}[theorem]{Problem}
\newtheorem{prop}[theorem]{Proposition}
\newtheorem{claim}{Claim}
\theoremstyle{remark}

\title{Almost regular subgraphs under spectral radius constrains}

\author{{Weilun Xu\thanks{Email: 225410010@fzu.edu.cn}}, {Guorong Gao\thanks{Email: grgao@fzu.edu.cn}}, {An Chang\thanks{Email: anchang@fzu.edu.cn}}
\\
\small{Center for Discrete Mathematics and Theoretical Computer Science\\ Fuzhou University\\ Fuzhou
Fujian, China}
 }

\date{ }

\begin{document}
\maketitle
\begin{center}
{\bf ABSTRACT}
\end{center}

A graph is called $K$-almost regular if its maximum degree is at most $K$ times the minimum degree. Erd\H{o}s and Simonovits showed that for a constant $0< \varepsilon< 1$ and a sufficiently large integer $n$, any $n$-vertex graph with more than $n^{1+\varepsilon}$ edges has a $K$-almost regular subgraph with  $n'\geq n^{\varepsilon\frac{1-\varepsilon}{1+\varepsilon}}$ vertices and at least $\frac{2}{5}n'^{1+\varepsilon}$ edges.
An interesting and natural problem is whether there exits the spectral counterpart to Erd\H{o}s and Simonovits's result. In this paper, we will completely settle this issue. More precisely, we verify that for  constants  $\frac{1}{2}<\varepsilon\leq 1$ and $c>0$, if the spectral radius of an $n$-vertex graph $G$ is at least $cn^{\varepsilon}$, then $G$ has a $K$-almost regular subgraph of order $n'\geq n^{\frac{2\varepsilon^2-\varepsilon}{24}}$ with at least $ c'n'^{1+\varepsilon}$ edges, where $c'$ and $K$ are constants depending on $c$ and $\varepsilon$. Moreover, for $0<\varepsilon\leq\frac{1}{2}$, there exist $n$-vertex graphs with spectral radius at least $cn^{\varepsilon}$ that do not contain such an almost regular subgraph.
Our result has a wide range of applications in spectral Tur\'{a}n-type problems. Specifically, let $ex(n,\mathcal{H})$ and $spex(n,\mathcal{H})$ denote, respectively, the maximum number of edges and  the maximum spectral radius among all $n$-vertex $\mathcal{H}$-free graphs. We show that for $1\geq\xi > \frac{1}{2}$, $ex(n,\mathcal{H}) = O(n^{1+\xi})$ if and only if $spex(n,\mathcal{H}) = O(n^\xi)$.

\section{Introduction}

All graphs considered in this paper are simple, i.e., graphs without parallel edges or loops.
For a given graph $G=(V,E)$  where $V=\{1,2,...,n\}$, the adjacency matrix $A(G)$ is defined as:
\[
A_{ij}(G) = \left\{
\begin{array}{ll}
    1 & \text{if } ij\in E(G), \\
    0  & \text{otherwise}.
\end{array}
\right.
\]
The spectral radius of $G$, which is denoted by $\lambda(G)$, corresponds to the largest eigenvalue of $A(G)$. Let $\Delta(G)$, $d(G)$ and $\delta(G)$ be the maximum degree, average degree and minimum degree of $G$, respectively.  A fundamental relationship among these parameters is well-established:
\begin{prop}
For any graph $G$, $$\delta(G) \leq d(G) \leq \lambda(G)\leq \Delta(G),$$ with equality holding in either case if and only if  $G$ is regular.
\end{prop}

Naturally, the gap between $\delta(G)$ and $\Delta(G)$ can serve as a measure of the "regularity" of the graph $G$.
In 1970,  Erd\H{o}s and Simonovits \cite{Erdos70} proposed the following definition of {\emph{$K$-almost regular}}, which occupies a pivotal role in the exploration of extremal graph theory.
\begin{dfn}[$K$-almost regular\cite{Erdos70}]
A graph $G$ is called $K$-almost regular for a constant $K\geq 1$ if $\Delta(G)\leq K\delta(G)$.
\end{dfn}


In \cite{Erdos70}, Erd\H{o}s and Simonovits proved that for sufficiently large $n$ and fixed $0<\varepsilon<1$, if an $n$-vertex graph has more than $n^{1+\varepsilon}$ edges, then $G$ contains an $n'$-vertex $K_\varepsilon$-almost regular subgraph with more than $\frac{2}{5}n'^{1+\varepsilon}$ edges, where $n'\geq n^{\varepsilon\frac{1-\varepsilon}{1+\varepsilon}}$, and $K_\varepsilon$ is a constant determined by $\varepsilon$. In \cite{Jiang}, Jiang and Seiver provided a variant of the result of Erd\H{o}s and Simonovits.

\begin{theorem}[Jiang and Seiver \cite{Jiang}, Proposition 2.7]\label{obs-1}
Let $\varepsilon$, $c$ be positive reals, with $\varepsilon \leq 1$. Let $n$ be a positive integer that is sufficiently large as a function of $\varepsilon$ and $c$. If $G$ is an $n$-vertex graph with $e(G)\geq cn^{1+\varepsilon}$, then $G$ contains a $K$-almost regular subgraph $G'$ on $n'\geq n^{\frac{\varepsilon}{4}}$ vertices such that $e( G^{\prime }) \geq c'n'^{1+ \varepsilon }$, where $K$ and $c'$ are constants depending on $c$ and $\epsilon$.
\end{theorem}
The results of \cite{Erdos70} and \cite{Jiang} reveal a fact that for sufficiently large $n$ and an $n$-vertex graph $G$, large average degree implies the existence of a large $K$-almost regular subgraph with large average degree. Inspired by this, we replace the average degree condition with a spectral condition and present the following theorem.
\begin{theorem}\label{Thm-1}
Let $G$ be a graph with order $n$ and $0<\varepsilon\leq \frac{1}{2}$, $c>0$ be two constants. Suppose that $n$ is sufficiently large and $\lambda(G)\geq cn^{\frac{1}{2}+\varepsilon}$. Then there exists a subgraph $G'$ of $G$ with order $n'\geq n^{\frac{\varepsilon}{3}}$ and a constant $c'=c(c,\varepsilon)$ such that $d(G')\geq c'n'^{\frac{1}{2}+\varepsilon}$.
\end{theorem}

Combining Theorem \ref{obs-1} and Theorem \ref{Thm-1}, we obtain the following theorem.
\begin{theorem}\label{Th-almost}
Let $G$ be a graph with order $n$ and $0<\varepsilon\leq \frac{1}{2}$, $c>0$ be two constants. Suppose that $n$ is sufficiently large and $\lambda(G)\geq cn^{\frac{1}{2}+\varepsilon}$. Then there exists a subgraph $G'$ of $G$ with order $n'\geq n^\frac{2\varepsilon^2+\varepsilon}{24}$ such that $G'$ is $K$-almost regular and $e(G')\geq c'n'^{\frac{3}{2}+\varepsilon}$, where $K$ and $c'$ are constants determined by $\varepsilon$ and $c$.
\end{theorem}

It is important to note that the bound of $\varepsilon>0$ in both Theorem \ref{Thm-1} and Theorem \ref{Th-almost} cannot be reduced. To emphasize this point, we present the following construction for any $0<\xi\leq \frac{1}{2}$.
\vskip 0.2 cm
\noindent\textbf{Construction:} Let $S_n^\xi$ be the graph obtained from the star graph $K_{1,\lceil n^{2\xi}\rceil-1}$ by subdividing one of its edges $n-\lceil n^{2\xi}\rceil$ times.

Clearly,  $S_n^\xi$ is an $n$-vertex graph. Note that  $S_n^\xi$ contains $K_{1,\lceil n^{2\xi}\rceil-1}$ as a subgraph. Thus, for $n$ large enough, $\lambda(S_n^\xi)\geq \lambda(K_{1,\lceil n^{2\xi}\rceil-1}) \geq\frac{1}{2}n^\xi$. However, any subgraph of $S_n^\xi$ is a linear forest. Hence, any subgraph of $S_n^\xi$ with order $n'\geq n^\frac{2\xi^2+\xi}{24}$ contains at most $n'-1$ edges.


\vskip 0.3 cm
Analogous to the widespread applications of the results of  Erd\H{o}s and Simonovits in classical Tur\'{a}n problems, our result likewise demonstrates broad utility in spectral Tur\'{a}n-type problems.

For a given graph class $\mathcal{H}$, a graph $G$ is called $\mathcal{H}$-free if it does not contain any graph from the class $\mathcal{H}$ as a subgraph. The \emph{ Tur\'{a}n number} (\emph{spectral Tur\'{a}n number}) of $\mathcal{H}$, denoted by $ex(n,\mathcal{H})$ ($spex(n,\mathcal{H})$), is equal to the maximum  number of edges (or maximum spectral radius for spectral Tur\'{a}n number) of an $n$-vertex $\mathcal{H}$-free graph. In particular, we write $ex(n,H)$ ($spex(n,H)$) instead of $ex(n,\mathcal{H})$ ($spex(n,\mathcal{H})$) when $\mathcal{H}=\{H\}$.

The classical Tur\'{a}n problem asks the value of $ex(n,\mathcal{H})$ and the extremal graph. In 1907 \cite{Mantel}, Mantel established that $ex(n,K_3)=\lfloor\frac{n^2}{4}\rfloor$, with the extremal graph being the complete bipartite graph $K_{\lfloor\frac{n}{2}\rfloor,\lceil\frac{n}{2}\rceil}$. This seminal work paved the way for further investigations. Tur\'{a}n \cite{Turan} extended Mantel's result to $K_{r+1}$-free graph, proving that the Tur\'{a}n graph, a complete $r$-partite graph with parts of equal or nearly equal sizes, maximizes the number of edges among all $n$-vertex $K_{r+1}$-free graphs. Later, Erd\H{o}s, Stone and Simonovits \cite{Erdos46}, \cite{Erdos66} gave a more general understanding of the Tur\'{a}n number for arbitrary graphs, which is known as the ESS Theorem. Specifically, for a given graph $H$ with chromatic number $\chi(H)\geq 2$, the ESS Theorem shows that
 $$ex(n,H)=(1-\frac{1}{\chi(H)-1}+o(1))\frac{n^2}{2}.$$
 ESS Theorem gives an approximate answer to the Tur\'{a}n number of a non-bipartite graph. However, for a bipartite graph $H$, the ESS Theorem implies only $ex(n,H)=o(n^2)$, leaving room for further exploration. Determining the order of magnitude of $ex(n,H)$ for some bipartite graph $H$ constitutes the \emph{degenerate Tur\'{a}n problem}, which is a prominent open problem in classical Tur\'{a}n problems.

A notable achievement in degenerate Tur\'{a}n problems was given by K\H{o}v\'{a}ri, S\'{o}s and Tur\'{a}n \cite{Kovari} in 1954, which bounded $ex(n,K_{s,t})$ as $O(n^{2-\frac{1}{s}})$ for $s\leq t$. Alon, Krivelevich and Sudakov \cite{Alon} generalized the result of K\H{o}v\'{a}ri, S\'{o}s and Tur\'{a}n to bipartite graphs satisfying all degrees in one partition class are at most $r$.
\begin{theorem}[Alon, Krivelevich and Sudakov \cite{Alon}]\label{thm-alon}
Let $H$ be a bipartite graph where all degrees in one partition class are at most $r$. Then $ex(n,H)=O(n^{2-\frac{1}{r}})$.
\end{theorem}
For a comprehensive overview of degenerate Tur\'{a}n problems, we recommend the survey by F\"{u}redi \cite{Furedi}.

The study of classical Tur\'{a}n problems naturally leads us to consider a related problem known as the \emph{spectral Tur\'{a}n-type problem}, which asks: What is the maximum spectral radius among all $n$-vertex $\mathcal{H}$-free graphs?

A seminal result in this area was given by Nikiforov \cite{Nikiforov07}, which established that $spex(n,K_{r+1})=\lambda(T_r(n))$, where $T_r(n)$ denotes the $r$-partite Tur\'{a}n graph on $n$ vertices. This result improved the result of Tur\'{a}n \cite{Turan}, because $d(G)\leq\lambda(G)$ holds for any graph $G$. Moreover, Nikiforov \cite{Nikiforov09} further developed a spectral analogue of ESS Theorem. Thus, the degenerate case also became  a major open problem on the research of spectral Tur\'{a}n-type problems. We refer readers to the classical survey \cite{Li} for more about spectral Tur\'{a}n-type problems.

It is natural to consider the relationship between $ex(n,\mathcal{H})$ and $spex(n,\mathcal{H})$. Notably, for any graph $G$, $d(G)\leq \lambda(G)$, which implies that if $spex(n,\mathcal{H})=O(n^c)$, then $ex(n,\mathcal{H})=O(n^{1+c})$. Conversely, if $ex(n,\mathcal{H})=\Omega(n^{1+c})$, then $spex(n,\mathcal{H})=\Omega(n^c)$.

Focusing on the extremal graphs of non-bipartite graphs, Wang, Kang and Xue \cite{Wang} proved that for sufficiently large $n$, if the graphs in $Ex(n,H)$ are Tur\'{a}n graphs plus $O(1)$ edges, then $Ex_{sp}(n,H)\subset Ex(n,H)$, where $Ex(n,H)=\{G:e(G)=ex(n,H), G~\text{is an}~n\text{-vertex}~H\text{-free}~\text{graph}\}$ and $Ex_{sp}(n,H)=\{G:\lambda(G)=spex(n,H), G ~\text{is an}~n\text{-vertex}~H\text{-free}~graph\}$.

In this paper, we establish a characterization of the relationship between $ex(n,\mathcal{H})$ and $spex(n,\mathcal{H})$ for more general $\mathcal{H}$.
\begin{theorem}\label{Thm-2}
Let $\frac{1}{2}<c\leq1$ be a constant. Then $ex(n,\mathcal{H})=O(n^{1+c})$ if and only if $spex(n,\mathcal{H})=O(n^c)$. Moreover, if $ex(n,\mathcal{H})=\Theta(n^{1+c})$, then $spex(n,\mathcal{H})=\Theta(n^c)$.
\end{theorem}

To the best of our knowledge, a spectral version of Theorem \ref{thm-alon} remains elusive. Consequently, for $r\geq 3$, we establish a spectral version of Theorem \ref{thm-alon}.
Clearly, Theorems \ref{thm-alon} and \ref{Thm-2} imply the following result.
\begin{theorem}
Let $H$ be a bipartite graph where all degrees in one partition class are at most $r\geq 3$. Then $spex(n,H)=O(n^{1-\frac{1}{r}})$.
\end{theorem}
Over the past decades, the exponents of  $ex(n,\mathcal{H})$ garnered significant attention. In \cite{Erdos81}, Erd\H{o}s conjectured that for any rational number $r\in [1,2]$, there exists a finite family of graphs $\mathcal{H}$ such that $ex(n,\mathcal{H})=(c_\mathcal{H}+o(1))n^{r}$ for some real number $c_\mathcal{H}$. In 2018, Bukh and Conlon \cite{Bukh} proved a weaker version of this conjecture.
\begin{theorem}[Bukh and Conlon\cite{Bukh}]\label{th-bc}
For every rational number $r\in [1,2]$, there exists a finite family of graphs $\mathcal{H}_r$ such that $ex(n,\mathcal{H}_r)=\Theta(n^r)$.
\end{theorem}
When restricting the finite family of graphs $\mathcal{H}$ to contain only one graph $H$, Erd\H{o}s \cite{Erdos88} further conjectured that for any rational number $r\in [1,2]$, there exists a graph $H$ satisfying $ex(n,H)=(c_H+o(1))n^{r}$ for some real number $c_H$. Kang, Kim and Liu\cite{Kang} subsequently stated a slightly weaker version of  Erd\H{o}s's conjecture.
\begin{conj}[Erd\H{o}s \cite{Erdos88}, Kang, Kim and Liu \cite{Kang}]\label{conj-liu}
For every rational number $r\in[1,2]$, there exists a graph $H$ such that $ex(n, H) = \Theta(n^r)$.
\end{conj}
Focusing on the exponents of $spex(n,\mathcal{H})$, we provide a spectral version of Theorem \ref{th-bc}.
\begin{theorem}\label{Thm-3}
For every rational number $r\in[\frac{1}{2},1]$, there exists a finite family of graphs $\mathcal{H}_r$ such that $spex(n,\mathcal{H}_r)=\Theta(n^r)$. Moreover, for every real number $0<r<\frac{1}{2}$, there is no finite family of graphs $\mathcal{H}$ such that $spex(n,\mathcal{H})=\Theta(n^r)$.
\end{theorem}
Moreover, if Conjecture \ref{conj-liu} is true, then it is easy to generalize to a spectral version.
\begin{theorem}
 If Conjecture \ref{conj-liu} is true, then for every rational number $r\in[\frac{1}{2},1]$, there exists a graph $H_r$ such that $spex(n,H_r)=\Theta(n^r)$. Moreover, for every rational number $0<r<\frac{1}{2}$, there is no graph $H$ such that $spex(n,H)=\Theta(n^r)$.
\end{theorem}
The rest of this paper is organized as follows. In Section 2, we introduce the notations and terminologies, and present the
proofs of our main results, following the order of Theorems \ref{Thm-1}, \ref{Thm-2} and \ref{Thm-3}. In Section 3, we make some conclusion remarks and propose a problem  that can be further researched.

\section{Proofs of main results}
Prior to embarking on the proofs, we introduce some notations and conclusions that will be frequently utilized in the subsequent proofs.

For a given $n$-vertex graph $G=(V,E)$, we assume $V=\{1,2,...,n\}$ without loss of generality. Let $A$ and $B$ be two subsets of $V(G)$. We denote by $E(A,B)$  the set of edges with one endpoint in $A$ and the other in $B$, and let $e(A,B)$ be its cardinality. When $A=B$, we simplify the notation to $E(A)$ and $e(A)$, respectively. Additionally, $G[A]$ represents the subgraph of $G$ induced by $A$, formally defined as $G[A]=(A,E(A))$.

Given an $n$-dimensional real vector $\mathbf{x}=\{\mathbf{x}_1,\mathbf{x}_2,...,\mathbf{x}_n\}^{\intercal}$, the Euclidean norm of $\mathbf{x}$ is defined as
$$\|\mathbf{x}\|_2=\sqrt{\sum_{i=1}^{n}\mathbf{x}_i^2}.$$
The celebrated Courant-Fischer Theorem asserts that the spectral radius of $G$, which is denoted by $\lambda(G)$, can be expressed as
$$\lambda(G)=\max_{\|\mathbf{x}\|_2=1}\mathbf{x}^{\intercal}A(G)\mathbf{x}=\max_{\|\mathbf{x}\|_2=1}2\sum_{uv\in E(G)}\mathbf{x}_u\mathbf{x}_v.$$
Moreover, there exists a non-negative unit eigenvector $\mathbf{x}$ of $A(G)$, associated with $\lambda(G)$ such that the equality $\lambda(G)=2\sum_{uv\in E(G)}\mathbf{x}_u\mathbf{x}_v$ holds.
\subsection{\bf Proof of Theorem \ref{Thm-1}.}\label{sec2.2}

 Let $G$ be an $n$-vertex graph with spectral radius $\lambda=\lambda(G)\geq cn^{\frac{1}{2}+\varepsilon}$, and let $\mathbf{x}$ be the non-negative eigenvector of $A(G)$ corresponding to $\lambda$ with $\|\mathbf{x}\|_2=1$ such that
\begin{equation}\label{eq-0}
\lambda=2\sum_{uv\in E(G)}\mathbf{x}_u\mathbf{x}_v.
\end{equation}
 Set $p=\lceil18^{\frac{2-2\varepsilon}{\varepsilon}}\rceil$ and $\eta=\frac{1}{p^2}$. We partition $V(G)$ into $p$ almost equal parts $B_1$,$B_2$,...,$B_p$, where $B_1$ contains $\lceil \frac{n}{p}\rceil$ vertices with the highest entries in $\mathbf{x}$.

We classify $G$ as being of $type$ $1$ if $\sum_{v\in B_1}\mathbf{x}_v^2\leq \frac{1}{2}-\eta$, and of $type$ $2$ otherwise.
\begin{claim}
    If $G$ is of type 1, then there exists a subgraph $G'$ of $G$ on $n'\geq \frac{p-1}{p}n-1$ vertices such that $d(G')\geq\frac{4c}{p^3}n'^{\frac{1}{2}+\varepsilon}$.
\end{claim}
\proof
Let $E_1=\{e\in E(G):e\cap B_1\neq \emptyset\}$. If $G$ is of type 1, then
\begin{equation}\label{eq-1}
(\frac{1}{2}-\eta)\lambda\geq\sum_{v\in B_1}\lambda\mathbf{x}_v^2= \sum_{v\in B_1}\mathbf{x}_v\sum_{u\sim v}\mathbf{x}_u\geq \sum_{uv\in E_1}\mathbf{x}_u\mathbf{x}_v.
\end{equation}
Consequently,
\begin{equation}\label{eq-2}
\sum_{uv\in E(G\setminus B_1)}\mathbf{x}_u\mathbf{x}_v= \sum_{uv\in E(G)}\mathbf{x}_u\mathbf{x}_v-\sum_{uv\in E_1}\mathbf{x}_u\mathbf{x}_v\geq \frac{1}{2}\lambda-(\frac{1}{2}-\eta)\lambda\geq \eta\lambda,
\end{equation}
where the second inequality is given by Inequalities (\ref{eq-0}) and (\ref{eq-1}).

Furthermore, by the definition of $B_1$, for every $v\in V(G)\setminus B_1$, we have $\mathbf{x}_v^2\leq \frac{\frac{1}{2}-\eta}{\frac{n}{p}} \leq \frac{p}{2}n^{-1}$. Therefore, for $uv\in E(G\setminus B_1)$,
\begin{equation}\label{eq-3}
\mathbf{x}_u\mathbf{x}_v\leq \sqrt{\frac{p}{2}n^{-1}}\times \sqrt{\frac{p}{2}n^{-1}}=\frac{p}{2}n^{-1}.
\end{equation}
Thus we have
\begin{equation}\notag
e(G\setminus B_1)\frac{p}{2}n^{-1}\geq \sum_{uv\in E(G\setminus B_1)}\mathbf{x}_u\mathbf{x}_v\geq \eta\lambda\geq \eta c n^{\frac{1}{2}+\varepsilon},
\end{equation}
where the first inequality is given by Inequality (\ref{eq-3}), and the second inequality is given by Inequality (\ref{eq-2}).

Hence, $e(G\setminus B_1)\geq \frac{2c\eta}{p} n^{\frac{3}{2}+\varepsilon}$. Set $G'=G\setminus B_1$ and $c'=\frac{4c\eta}{p}=\frac{4c}{p^3}$. We have $n'=|V(G')|\geq \lfloor\frac{p-1}{p}n\rfloor$ and $d(G')\geq c'n'^{\frac{3}{2}+\varepsilon}$.
\endproof
Consequently, if $G$ is of type 1, then the proof is completed. Thus we assume that $G$ is of type 2.
\begin{claim}
    If $G$ is of type 2, then there exists a subgraph $G_1$ of $G$, such that $\frac{n}{p}\leq |V(G_1)|\leq \frac{3n}{p}$ and $\lambda(G_1)\geq \frac{\lambda}{6\sqrt p}\geq c|V(G_1)|^{\frac{1}{2}+\varepsilon}$.
\end{claim}
\proof
Since $G$ is of type 2, we have $\sum_{v\in B_1}\mathbf{x}_v^2> \frac{1}{2}-\eta$. Therefore,
\begin{align}\notag
\sum_{v\in B_1}\lambda \mathbf{x}_v^2 &= \sum_{v\in B_1}\mathbf{x}_v\sum_{u\sim v}\mathbf{x}_u\\ \notag
&=2\sum_{uv\in E(G[B_1])}\mathbf{x}_u\mathbf{x}_v+\sum_{i= 2}^p(\sum_{uv\in E(G(B_1,B_i))}\mathbf{x}_u\mathbf{x}_v)\\ \notag
&> (\frac{1}{2}-\eta)\lambda.\notag
\end{align}
Let $$f=\sum_{uv\in E(G[B_1])}\mathbf{x}_u\mathbf{x}_v.$$
Then we have
\begin{equation}\label{eq-gamma}
\sum_{i=2}^{p}\sum_{uv\in E(G(B_1,B_i))}\mathbf{x}_u\mathbf{x}_v\geq (\frac{1}{2}-\eta)\lambda-2f.
\end{equation}
For $2\leq i\leq p$, suppose that
$$\alpha_i=\sum_{v\in V(G)\setminus(B_1\cup B_i)}\mathbf{x}_v^2,~\beta_i=\sum_{v\in B_i}\mathbf{x}_v^2 ~~\text{and}~~\gamma_i=\sum_{uv\in E(G(B_1,B_i))}\frac{\mathbf{x}_u\mathbf{x}_v}{\lambda}.$$

\textbf{Case 1:} $f\geq \frac{\lambda}{\sqrt {p}}$.

 In this case, we simply  define $G_1$ as $G[B_1]$. Then $|V(G_1)|=\lceil\frac{n}{p}\rceil\leq \frac{n}{p}+1\leq \frac{2n}{p}$ and $\lambda(G_1)\geq \frac{\lambda}{\sqrt {p}}$.

\textbf{Case 2:} $f< \frac{\lambda}{\sqrt {p}}$, and there exists an integer $2\leq j\leq p$ such that $\alpha_j\leq \frac{1}{2}-\frac{1}{\sqrt p}$.

 Let $G_1=G[B_1\cup B_j]$ and $E_2=E(G)\setminus E(G_1)$. Then any $e\in E_2$ contains at least one end point in $V(G)\setminus(B_1\cup B_j)$. Thus, similar to Inequality (\ref{eq-1}), we have
\begin{align}\notag
(\frac{1}{2}-\frac{1}{\sqrt p})\lambda \geq \alpha_j\lambda = \sum_{v\in V(G)\setminus (B_1\cup B_j)}\lambda\mathbf{x}_v^2=\sum_{v\in V(G)\setminus (B_1\cup B_j)}\mathbf{x}_v\sum_{u\sim v}\mathbf{x}_u\geq \sum_{uv\in E_2}\mathbf{x}_u\mathbf{x}_v.
\end{align}
Consequently, similar to Inequality (\ref{eq-2}) we have
\begin{align*}
\sum_{uv\in E(G_1)}\mathbf{x}_u\mathbf{x}_v=\sum_{uv\in E(G)}\mathbf{x}_u\mathbf{x}_v-\sum_{uv\in E_2}\mathbf{x}_u\mathbf{x}_v\geq \frac{\lambda}{2}-(\frac{1}{2}-\frac{1}{\sqrt p})\lambda=\frac{\lambda}{\sqrt p}.
\end{align*}
Clearly, $|V(G_1)|\leq 2\lceil\frac{n}{p}\rceil\leq \frac{2n}{p}+2\leq \frac{3n}{p}$. Moreover, by Courant-Fischer Theorem, we have $\lambda(G_1)\geq \frac{\lambda}{\sqrt{p}}$.

\textbf{Case 3:} $f< \frac{\lambda}{\sqrt {p}}$ and $\alpha_i\geq \frac{1}{2}-\frac{1}{\sqrt p}$ for any $2\leq i\leq p$.

In this case, by Inequality (\ref{eq-gamma}) and the definition of $\gamma_i$, we have
\begin{equation}\label{eq-6}
\sum_{i=2}^{p}\gamma_i\geq (\frac{1}{2}-\eta-\frac{2}{\sqrt p})\geq \frac{1}{3}.
\end{equation}
Moreover, since $G$ is of type 2,
\begin{equation}\label{eq-5}
\sum_{i=2}^{p}\beta_i=1-\sum_{v\in B_1}\mathbf{x}_v^2\leq 1-(\frac{1}{2}-\eta)\leq \frac{2}{3}.
\end{equation}

We first show that there exists an integer $2\leq j\leq p$ such that
\begin{equation}\label{eq-7}
\frac{\gamma_j}{\sqrt{3\beta_j}}> \frac{1}{9\sqrt{p}}.
\end{equation}
Suppose to the contrary that $\frac{\gamma_i}{\sqrt{3\beta_i}}\leq \frac{1}{9\sqrt{p}}$ for every $2\leq i\leq p$. Then for any $2\leq i\leq p$, $9\sqrt{p}\gamma_i\leq \sqrt{3\beta_i}$. Thus,
$$ 3\sqrt{p}\leq \sum_{i=2}^{p}9\sqrt{p}\gamma_i\leq \sum_{i=2}^{p}\sqrt{3\beta_i}\leq \sqrt{2(p-1)},$$
where the first inequality is given by Inequality (\ref{eq-6}), and the last inequality is given by Inequality (\ref{eq-5}) and Cauchy-Schwarz Inequality. However, the definition of $p$ implies that $3\sqrt p>\sqrt{2(p-1)}$, which is a contradiction.
Let $j$ be the integer satisfying Inequality (\ref{eq-7}), and $G_1=G(B_1,B_j)$. Then $|V(G_1)|\leq 2\lceil \frac{n}{p}\rceil\leq\frac{3n}{p}$. Moreover,
\begin{equation}
\sum_{uv\in E(G_1)}\mathbf{x}_u\mathbf{x}_v=\gamma_j\lambda.
\end{equation}
Now, define a $|B_1\cup B_j|$-dimensional vector $\mathbf{y}$ as following:
\[
\mathbf{y}_v = \left\{
\begin{array}{ll}
    \mathbf{x}_v & \text{if } v\in B_1, \\
    \sqrt{\frac{\alpha_j+\beta_j}{\beta_j}}\mathbf{x}_v  & \text{if } v\in B_j.
\end{array}
\right.
\]
Clearly, $\|\mathbf{y}\|_2=1$. Recall that in this case, $\alpha_j\geq \frac{1}{2}-\frac{1}{\sqrt p}\geq \frac{1}{3}$. So we have $\mathbf{y}_v\geq \frac{1}{\sqrt{3\beta_j}}\mathbf{x}_v$ for any $v\in B_j$.
Moreover, every edge in $G_1$ contains an end point in $B_j$.
By Inequality (\ref{eq-7}) and Courant-Fischer Theorem, we have
\begin{equation}\notag
\lambda(G_1)\geq 2\sum_{uv\in E(G_1)}\mathbf{y}_u\mathbf{y}_v\geq \frac{2}{\sqrt{3\beta_j}}\sum_{uv\in E(G_1)}\mathbf{x}_u\mathbf{x}_v=\frac{2\gamma_j}{\sqrt{3\beta_j}}\lambda\geq\frac{2\lambda}{9\sqrt{p}}.
\end{equation}
 Consequently, if $G$ is of type 2, then there is a subgraph $G_1$ of $G$ such that $\frac{n}{p}\leq |V(G_1)|\leq\frac{3n}{p}$ and $\lambda(G_1)\geq\frac{2\lambda}{9\sqrt p}\geq\frac{\lambda}{6\sqrt p}$. Moreover, by the definition of $p$, $p^{\varepsilon}>18$. Hence,
 $$\lambda(G_1)\geq \frac{c}{6\sqrt p}n^{\frac{1}{2}+\varepsilon}\geq \frac{c p^{\varepsilon}}{6\times 3^{\frac{1}{2}+\varepsilon}}|V(G_1)|^{\frac{1}{2}+\varepsilon} \geq c|V(G_1)|^{\frac{1}{2}+\varepsilon}.$$
\endproof
Now, we can replace $G$ with $G_1$ and repeat the above process. If $G_1$ is of type 2, we define $G_2$ from $G_1$ the way we define $G_1$ from $G$. We stop the process once $G_i$ is of type 1. Suppose that $k$ is the maximum integer such that $G_k$ is of type 2, then
$$|V(G_k)|\leq (\frac{3}{p})^kn,~\text{and} ~\lambda(G_k)\geq (\frac{1}{6\sqrt p})^k\lambda\geq (\frac{1}{6\sqrt p})^kcn^{\frac{1}{2}+\varepsilon}.$$
Note that for any graph $H$, $\lambda(H)< |V(H)|$. Thus we have
$$(\frac{1}{6\sqrt p})^kcn^{\frac{1}{2}+\varepsilon}\leq(\frac{3}{p})^kn.$$
 Then $(\frac{\sqrt p}{18})^k\leq \frac{1}{c}n^{\frac{1}{2}-\varepsilon}$, which implies that
$$k\leq ((1-2\varepsilon)\log n-\log c)/\log \frac{p}{18^2}.$$
Note that $|V(G_k)|\geq \frac{n}{p^k}$,
$$\log|V(G_k)|\geq (1-(1-2\varepsilon)\frac{\log p}{\log\frac{p}{18^2}}+\frac{\log c \log p}{\log\frac{p}{18^2}\log n})\log n.$$
 Since $n$ is sufficiently large, we have $\mid\frac{\log c \log p}{\log\frac{p}{18^2}\log n}\mid\leq \frac{\varepsilon}{2}$. Hence
  $$\log|V(G_k)|\geq (1-(1-2\varepsilon)\frac{\log p}{\log\frac{p}{18^2}}-\frac{\varepsilon}{2})\log n.$$
   By the definition of $p$, we have $|V(G_k)|\geq n^{\frac{\varepsilon}{2}}$. Then $G_{k+1}$ is of type 1 and $|V(G_{k+1})|\geq \frac{n^{\frac{\varepsilon}{2}}}{p}\geq n^{\frac{\varepsilon}{3}}$. Consequently, there is a subgraph $G'$ of $G_{k+1}$ with order $n'\geq\frac{p-1}{p}n^{\frac{\varepsilon}{3}}$, such that $d(G')\geq c'n'^{\frac{1}{2}+\varepsilon}$, where $c'=\frac{4c}{p^3}$.
\endproof
\subsection{\bf Proof of Theorem \ref{Thm-2}.}

We first show that $ex(n,\mathcal{H})=\Omega(n^{1+c})$ implies $spex(n,\mathcal{H})=\Omega(n^c)$. Since $ex(n,\mathcal{H})=\Omega(n^{1+c})$, there exist constants $C$, $N>0$, such that for all $n>N$, there exists an $\mathcal{H}$-free graph $G$ with $e(G)\geq Cn^{1+c}$. Hence, $\lambda(G)\geq d(G)= \frac{2e(G)}{n}\geq 2Cn^c$, implying $spex(n,\mathcal{H})=\Omega(n^c)$.

Next, we prove that $ex(n,\mathcal{H})=O(n^{1+c})$ implies $spex(n,\mathcal{H})=O(n^c)$. Suppose to the contrary that for any constant $C>0$, there exists a sufficiently large integer $n_1$ and an $n_1$-vertex $\mathcal{H}$-free graph $G$ with $\lambda(G)\geq Cn_1^{c}$. Then by Theorem \ref{Thm-1}, there exists a subgraph $G'$ of $G$, such that
 $$|V(G')|=n'\geq n_1^{\frac{c}{3}-\frac{1}{6}}, ~\text{and}~ d(G')\geq C'n'^{c}.$$
Since $ex(n,\mathcal{H})=O(n^{1+c})$,
there exist constants $C''$ and $N>0$, such that for all $n\geq N$, any $n$-vertex $\mathcal{H}$-free graph has at most $C''n^{1+c}$ edges. Choose $C$ and $n_1$ such that $C'>\frac{C''+1}{2}$ and $n_1^{\frac{c}{3}-\frac{1}{6}}\geq N$. Then there is an $\mathcal{H}$-free graph with order $n'\geq n_1^{\frac{c}{3}-\frac{1}{6}}\geq N$ and at least $(C''+1)n'^{1+c}>C''n'^{1+c}$ edges, contradicting $ex(n,\mathcal{H})=O(n^{1+c})$.
\endproof

\subsection{\bf Proof of Theorem \ref{Thm-3}.}

The first part of the theorem  directly follows from Theorems \ref{Thm-2} and \ref{th-bc}, as well as the known result in \cite{Desai} that $spex(n,C_{2k})=\Theta(n^{\frac{1}{2}})$. Therefore, we focus on proving that for any finite graph class $\mathcal{H}$, there does not exist a constant $c\in(0,\frac{1}{2})$ such that $spex(n,\mathcal{H})=\Theta(n^c)$.

We consider two cases:

\textbf{Case 1:} $\mathcal{H}$ contains $K_{1,p} \cup qK_1$ for some $p, q \in \mathbb{N}$. \\
In this case, for sufficiently large $n$, any $n$-vertex $\mathcal{H}$-free graph has maximum degree at most $p-1$. Hence, the spectral radius of such a graph is bounded by $p-1$, which implies $spex(n,\mathcal{H}) = O(1)$. Clearly, this rules out the existence of a constant $c \in (0,\frac{1}{2})$ such that $spex(n,\mathcal{H})=\Theta(n^c)$.

\textbf{Case 2:} $\mathcal{H}$ does not contain $K_{1,p} \cup qK_1$ for any $p, q \in \mathbb{N}$. \\
In this case, the complete bipartite graph $K_{1,n-1}$ is $\mathcal{H}$-free. Since the spectral radius of $K_{1,n-1}$ is $\sqrt{n-1}$, it follows that $spex(n,\mathcal{H}) \geq \sqrt{n-1} = \Omega(n^{\frac{1}{2}})$. Again, this rules out the existence of a constant $c \in (0,\frac{1}{2})$ such that $spex(n,\mathcal{H})=\Theta(n^c)$.

In either case, we have shown that for any finite graph class $\mathcal{H}$, there is no constant $c \in (0,\frac{1}{2})$ satisfying $spex(n,\mathcal{H})=\Theta(n^c)$.
\section{Concluding Remark}

The original version of Theorem \ref{obs-1} (Proposition 2.7 in \cite{Jiang}) requires $c\geq1$ and $\varepsilon<1$. However,
 the requirement of $c\geq 1$ and $\varepsilon<1$ is not essential. Indeed, in the proof of Proposition 2.7 \cite{Jiang}, a more sophisticated calculation shows the existence of a $K$-almost regular subgraph with $m\geq n^{\frac{\varepsilon}{4}}$ vertices, as long as $\mid\frac{\log c \log p}{\log \frac{p}{4} \log n}\mid\leq \frac{\varepsilon}{4}$ and $\frac{\log p}{\log \frac{p}{4}}\leq \frac{2-\varepsilon}{2-2\varepsilon}$ ($p\geq 2^{\frac{4-2\varepsilon}{\varepsilon}}$). This is similar to the last part of our proof of Theorem \ref{Thm-1} in Section \ref{sec2.2}.

The lower bound of the average degree of the subgraph $G'$ in Theorem \ref{Thm-1} cannot be improved to $c'n^{\frac{1}{2}+\varepsilon}$. In fact, the graph $K_{n^{2\varepsilon},n-n^{2\varepsilon}}$ has spectral radius $\sqrt{n^{2\varepsilon}(n-n^{2\varepsilon})}\geq\frac{1}{2}n^{\frac{1}{2}+\varepsilon}$. However, any subgraph of this graph has average degree at most $n^{2\varepsilon}$.

In this paper, for sufficiently large $n$, we proved that if an $n$-vertex graph $G$ satisfying $\lambda(G)\geq cn^{\frac{1}{2}+\varepsilon}$, then there is a sufficiently large almost regular subgraph $G'$ of $G$ with $\Omega(|V(G')|^{\frac{1}{2}+\varepsilon})$ edges. Conversely, we also showed that there exist an $n$-vertex graph $G$ with $\lambda(G)=\Theta(n^{\frac{1}{2}})$ that does not possess any such an almost regular subgraph. However, the case where $\lambda(G)$ is "slightly" larger than $n^{\frac{1}{2}}$, specifically when $\lambda(G)\geq n^{\frac{1}{2}}\log n$, remains unclear. This prompts the following problem:
\begin{problem}
For sufficiently large $n$ and an $n$-vertex graph $G$ with $\lambda(G)\geq n^{\frac{1}{2}}\log n$, whether there is a constant $K$ independent of $n$ such that $G$ contains a sufficiently large $K$-almost regular subgraph $G'$ with $\Omega(|V(G')|^{\frac{1}{2}}\log|V(G')|)$ edges?
\end{problem}

\end{document}